\documentclass[10pt]{article}

\usepackage{amscd,amsmath, amssymb, fancyhdr, mathbbol}

\numberwithin{equation}{section}

\newcommand{\version}{version 2.4.1,\ \ November 13, 2018}

\def\eqref#1{(\ref{#1})}
\newcommand{\goth}{\mathfrak}

\newcommand{\arrow}{{\:\longrightarrow\:}}
\newcommand{\Z}{{\Bbb Z}}
\def\C{{\Bbb C}}
\newcommand{\R}{{\Bbb R}}

\newcommand{\6}{\partial}
\def\1{\sqrt{-1}\:}
\newcommand{\restrict}[1]{{\left|_{{\phantom{|}\!\!}_{#1}}\right.}}
\newcommand{\cntrct}                
{\hspace{2pt}\raisebox{1pt}{\text{$\lrcorner$}}\hspace{2pt}}

\renewcommand{\tilde}{\widetilde}
\renewcommand{\bar}{\overline}
\renewcommand{\phi}{\varphi}
\renewcommand{\epsilon}{\varepsilon}
\renewcommand{\geq}{\geqslant}
\renewcommand{\leq}{\leqslant}
\renewcommand{\max}{{\rm max}}

\newcommand{\Pic}{\operatorname{Pic}}
\newcommand{\Par}{\operatorname{Par}}

\newcommand{\Aut}{\operatorname{Aut}}
\newcommand{\Alb}{\operatorname{Alb}}

\newcounter{Mycounter}[section]
\newcounter{lemma}[section]
\renewcommand{\thelemma}{{Lemma \thesection.\arabic{lemma}}}
\newcommand{\lemma}{%
    \setcounter{lemma}{\value{Mycounter}}
    \refstepcounter{lemma}
    \stepcounter{Mycounter}
    {\noindent \bf \thelemma:\ }}

\newcounter{claim}[section]
\newcounter{sublemma}[section]
\newcounter{corollary}[section]
\newcounter{theorem}[section]
\renewcommand{\thetheorem}{{Theorem \thesection.\arabic{theorem}}}
\newcommand{\theorem}{%
    \setcounter{theorem}{\value{Mycounter}}
    \refstepcounter{theorem}
    \stepcounter{Mycounter}
    {\noindent \bf \thetheorem:\ }}

\newcounter{conjecture}[section]
\newcounter{proposition}[section]
\renewcommand{\theproposition}
      {{Proposition \thesection.\arabic{proposition}}}
\newcommand{\proposition}{%
    \setcounter{proposition}{\value{Mycounter}}
    \refstepcounter{proposition}
    \stepcounter{Mycounter}
    {\noindent \bf \theproposition:\ }}

\newcounter{definition}[section]
\renewcommand{\thedefinition}
      {{Definition~\thesection.\arabic{definition}}}
\newcommand{\definition}{%
    \setcounter{definition}{\value{Mycounter}}
    \refstepcounter{definition}
    \stepcounter{Mycounter}
    {\noindent \bf \thedefinition:\ }}

\newcounter{example}[section]
\newcounter{remark}[section]
\renewcommand{\theremark}{{Remark \thesection.\arabic{remark}}}
\newcommand{\remark}{%
    \setcounter{remark}{\value{Mycounter}}
    \refstepcounter{remark}
    \stepcounter{Mycounter}
    {\noindent \bf \theremark:\ }}

\newcounter{problem}[section]
\newcounter{question}[section]
\renewcommand{\thequestion}{{Question \thesection.\arabic{question}}}
\newcommand{\question}{%
    \setcounter{question}{\value{Mycounter}}
    \refstepcounter{question}
    \stepcounter{Mycounter}
    {\noindent \bf \thequestion:\ }}

\newcommand{\proof}{\noindent{\bf Proof:\ }}

\makeatletter

\@addtoreset{equation}{section} \@addtoreset{footnote}{section}
\makeatother

\def\blacksquare{\hbox{\vrule width 5pt height 5pt depth 0pt}}
\def\endproof{\blacksquare}

\begin{document}
\begin{center}
{\LARGE\bf
Algebraically hyperbolic manifolds have 
finite automorphism groups\\[4mm]
}

Fedor Bogomolov, Ljudmila Kamenova\footnote{Partially supported 
by a grant from the Simons Foundation/SFARI (522730, LK)}, 
Misha Verbitsky\footnote{Partially supported 
by the  Russian Academic Excellence Project '5-100', 
CNPq Process 313608/2017-2 and FAPERJ E-26/202.912/2018.
}

\end{center}

{\small \hspace{0.10\linewidth}
\begin{minipage}[t]{0.85\linewidth}
{\bf Abstract.} 
A projective manifold $M$ is algebraically hyperbolic if there exists a 
positive constant $A$ such that the degree of any curve of genus $g$ on $M$ 
is bounded from above by $A(g-1)$. A classical result is that Kobayashi 
hyperbolicity implies algebraic hyperbolicity. It is known that Kobayashi 
hyperbolic manifolds have finite automorphism groups. Here we prove that, 
more generally, algebraically hyperbolic projective manifolds have finite 
automorphism groups. 
\end{minipage}
}


\section{Introduction} 

J.-P. Demailly introduced an algebraic analogue for the analytic notions 
of hyperbolicity in \cite{Demailly}. Using the classical 
Gauss-Bonnet formula he showed that every compact complex Kobayashi 
hyperbolic manifold is also algebraically hyperbolic. Demailly conjectured 
that for complex projective manifolds the two notions, Kobayashi 
hyperbolicity and algebraic hyperbolicity, coincide. 

\hfill

Kobayashi proved that Kobayashi hyperbolic manifolds have only finite 
order automorphisms, \cite{_Kobayashi:1976_}. In \cite{kv} L. Kamenova and 
M. Verbitsky proved that algebraically hyperbolic projective 
hyperk\"ahler manifolds have finite groups of automorphisms. The general
expectation is that all projective hyperk\"ahler manifolds are algebraically 
non-hyperbolic.

\hfill

\theorem (Kamenova, Verbitsky)
Let $M$ be an algebraically hyperbolic projective hyperk\"ahler manifold. 
Then $M$ has finite group of automorphisms. 

\hfill

Here we generalize the result above and show that all algebraically 
hyperbolic projective manifolds have finite groups of automorphisms. 

\hfill

\theorem 
The group $\Aut(M)$ of automorphisms of an algebraically hyperbolic 
manifold $M$ is finite. 

\hfill

The idea behind the proof of our main theorem is the following. 
We consider the induced action of 
the automorphism group $\Aut(M)$ on the cohomology $H^{1,1}(M,\R)$. 
From the Hodge-Riemann relations it follows that the image of $\Aut(M)$ 
in $GL(H^{1,1}(M,\R))$ has to preserve some rational K\"ahler class. 
Furthermore, we consider 
the case when the image of $\Aut(M)$ in $\Aut(\Pic^0(M))$ is infinite. 
Using the Albanese map we produce a subvariety of $M$ that admits 
self-isogenies of arbitrarily high order, thus giving curves of constant 
genus and arbitrary high degree in this subvariety, i.e., we obtain an 
algebraically non-hyperbolic subvariety. This contradiction implies 
that the image of $\Aut(M)$ in $\Aut(\Pic^0(M))$ is finite. If we assume 
that $\Aut(M)$ is infinite, there would be an infinite subgroup $\Gamma$ 
acting trivially on $\Aut(\Pic(M))$. It fixes an ample line bundle 
$L \in \Pic(M)$. The action of $\Gamma$ on ${\Bbb P} H^0(M,L)^*$ preserves the 
image of the projective embedding $M \hookrightarrow {\Bbb P} H^0(M,L)^*$. 
Let $\bar\Gamma$ be the Zariski closure of $\Gamma$ in $PGL(H^0(M,L)^*)$. 
Then the orbits of $\bar\Gamma$ are positive dimensional, which is a 
contradiction, because $\Aut(M)$ is discrete when $M$ is algebraically 
hyperbolic. 


\section{Brody curves, Kobayashi and algebraic hyperbolicity}


In this section let $M$ be a compact complex manifold. 
We introduce some basic hyperbolicity notions that can be found in 
\cite{_Kobayashi:1976_} and \cite{Demailly}. 
Brody introduced the Brody curves in \cite{_Brody:hyperboloc_}.

\subsection{Kobayashi pseudometric and Brody curves}

\hfill

\definition
A {\it pseudometric} on $M$ is a function $d:\; M \times M \arrow \R^{\geq 0}$ 
which is symmetric: $d(x,y)=d(y,x)$ and satisfies the 
triangle inequality $d(x,y)+d(y,z) \geq d(x,z)$. 

\hfill 

\remark  
Let ${\goth D}$ be a set of pseudometrics on $M$.  Then 
$d_{\max}(x,y):= \sup_{d\in {\goth D}}d(x,y)$ is also a pseudometric. 


\hfill

\definition 
The {\it Kobayashi pseudometric} $d_M$ on a complex manifold $M$ 
is $d_M:=d_{\max}$ for the set ${\goth D}$ of all pseudometrics on $M$ 
such that any holomorphic map from the 
Poincar\'e disk to $M$ is distance-decreasing. 

\hfill

\remark 
The definition above is equivalent to the standard definition that 
the Kobayashi distance between points $x, y \in M$ 
is the infimum of the Poincar\'e distance over all sets of Poincar\'e disks
connecting $x$ to $y$. 

\hfill







\definition
A manifold $M$ is called {\it Kobayashi hyperbolic} if 
the Kobayashi pseudometric $d_M$ is non-degenerate. 

\hfill 




\definition
Let $M$ be a complex Hermitian manifold. 
A {\it Brody curve} is a non-constant 
holomorphic map $f:\; \C\arrow M$ 
such that $|df|\leq C$ for some constant $C$. 
Here $|df|$ is understood as an operator norm of 
$df:\; T_z \C \arrow TM$, where $\C$ is equipped 
with the standard Euclidean metric. 

\hfill

\definition
Let $(\Delta_r, g_r)$ be a disk of radius $r$ in $\C$ 
with the Poincar\'e metric $g_r$, rescaled in such a way that 
the unit tangent vector to $0$ has length $1$. 
A {\it Brody map} to a Hermitian complex manifold $M$ 
is a map $f:\; \Delta_r\arrow M$ such that $|df| \leq 1$ 
(here the operator norm is taken with respect to the 
Poincar\'e metric on $\Delta_r$) and $|df|(z)=1$ at $z=0$. 

\hfill

\lemma \label{Lemma1} 
Let $f_r:\; \Delta_r\arrow M$ be a sequence of Brody maps 
with $r \arrow \infty$. Then $\{ f_r \}$ has a subsequence which 
converges uniformly to 
a Brody curve $f$ satisfying $|df|(z)=1$ at $z=0$. 

\hfill

{\bf Proof}. Let $r_1 < r_2$. The identity map 
$\tau:\; (\Delta_{r_1}, g_1) \arrow (\Delta_{r_2}, g_2)$
is 1-Lipschitz. Indeed, it is 1-Lipschitz 
with respect to the usual Poincar\'e metric:
$\tau^* (r_2^{-2}g_2) \leq r_1^{-2} g_1$. 
Since $r_1 < r_2$, this gives $\tau^* g_2 \leq g_1$.
Restricted on any disk $\Delta_R$, the family $\{f_r, r> R\}$ is a normal family
(since it is Lipschitz), hence it has a subsequence which converges uniformly to a Lipschitz map.
Since a uniform limit of holomorphic maps is holomorphic,
the family $\left\{f_r\restrict{\Delta_R}, r> R\right\}$ converges to 
a holomorhic map on $\phi_R:\; \Delta_R\arrow M$. 
The map $\phi_R$ is Lipschitz with respect to all metrics $g_r$, $r>R$. 
Since $\lim\limits_{r \arrow \infty} g_r$ is the standard Euclidean metric 
$g_\infty$, the map $\phi_R$ is Lipschitz with respect to $g_\infty$. 
The sequence $\lim\limits_{R \arrow \infty} \phi_R$ converges to a 
holomorphic Lipschitz map $\C \arrow M$. Since all $f_r$ and $\phi_R$ satisfy 
$|d\phi_R|(z)=1$ at $z=0$, the same is true for the limit. 
\endproof

\hfill

A classical result of Brody states that a compact complex manifold $M$ 
is Kobayashi hyperbolic if and only if there are no non-constant 
holomorphic maps from $\C$ to $M$. 

\hfill 

\theorem(Brody's lemma) \label{brodyslemma}
Let $M$ be a compact complex manifold which is not 
Kobayashi hyperbolic. Then $M$ contains a Brody curve. 

\hfill

\proof
Let us equip $M$ with a Hermitian metric $h$. 
If $|df|(0)\leq C$ for any holomorphic 
map $(\Delta_1, g_1)\arrow M$, then the Kobayashi 
pseudometric satisfies $d_K\geq C^{-1}h$, and $M$ is Kobayashi 
hyperbolic. If this quantity is non-bounded, we can always rescale 
the disc to obtain a map $f_r:\; (\Delta_r, g_r)\arrow M$ 
with $r= |df|(0)$, and then $|df_r|(0)=1$.
Then Brody's lemma follows from \ref{Lemma1} and the following lemma. 

\hfill

\lemma
Let $M$ be a compact Hermitian manifold,
and $\psi_r:\; (\Delta_r, g_r)\arrow M$ a sequence
of holomorphic maps satisfying $|d\psi_r|(0)\geq 1$, $r\arrow \infty$.
Then there exists a sequence of Brody maps 
$f_s:\; (\Delta_s, g_s)\arrow M$, with $s\arrow \infty$.

\hfill

{\bf Proof}. We need to construct a sequence of Brody maps,
which are 1-Lipschitz maps $f_s:\; (\Delta_s, g_s)\arrow M$, with
$|df_s|(0)=1$. The identity map 
\[ \Psi_{r-\epsilon, r}:\; (\Delta_{r-\epsilon},g_{r-\epsilon})\arrow (\Delta_r, g_r)
\]
is 1-Lipschitz, and satisfies 
\[ \lim_{z\arrow \6\Delta_{r-\epsilon}}|d\Psi_{r-\epsilon,
  r}|(z)=0.
\]
Let $u:= r-\epsilon$ and $\tilde f_u:= \Psi_{r-\epsilon,
  r}\circ\psi_r$ be a restriction of $f_r$ to the disk
$(\Delta_{r-\epsilon},g_{r-\epsilon})$.
Then $f_u$ is also Lipschitz and 
$|d\tilde f_u|$ reaches maximum at a point $z_u$
somewhere inside the disk $\Delta_u$.
Applying an appropriate holomorphic
isometry of $\Delta_u$, we may assume that $|d\tilde f_u|(z)$
takes maximum $C_u\geq 1$ for $z=0$. Rescaling $\tilde f_u$, and
setting $s:= C_uu$, we obtain a map $f_s:\; \Delta_s\arrow M$ 
which is 1-Lipschitz and satisfies $|df_s|\leq 1$,
$|df_s|(0)=1$.
\endproof

\hfill

\definition
Let $M$ be a projective manifold. 
We say that $M$ is {\it algebraically hyperbolic} 
if there exists a constant $A >0$ such that for any curve 
$C\subset M$ of geometric genus $g$ one has $\deg C < A (g-1)$. 

\hfill

\remark Algebraically hyperbolic manifolds contain no 
elliptic nor rational curves.

\hfill

Using Gauss-Bonnet in \cite{Demailly} Demailly proved that Kobayashi 
hyperbolicity implies algebraic hyperbolicity. Here we'll give a 
slightly different proof of this well known result. The converse implication, 
``algebraically hyperbolicity implies Kobayashi hyperbolicity'', was 
conjectured by J.-P. Demailly who introduced the notion of algebraic 
hyperbolicity in \cite{Demailly}. 

\hfill 

\theorem (Demailly, \cite{Demailly}) 
Kobayashi hyperbolicity implies algebraic hyperbolicity. 

\hfill

{\bf Proof}. Let $C\subset M$ be a curve in a Kobayashi hyperbolic manifold. 
The genus $0$ and $1$ cases are ruled out by Brody's Lemma, \ref{brodyslemma}. 
Then genus $g=g(C)$ is greater than $1$, and its universal cover is the disk 
$\Delta$. Denote by $\phi_C:\; \Delta \arrow M$ the universal covering map. 
The volume of $C$ with respect to the Fubini-Study metric on $M$ is $\deg C$, 
and its volume with respect to the Poincar\'e metric is 
$\int_C c_1(T^*C)=2\pi(2g-2)$. Therefore, 
$|\phi_C| > \frac{\deg C}{2\pi(2g-2)}$ somewhere on $\Delta$. 
Notice that $M$ is not algebraically hyperbolic if and only if 
there is a sequence $\{ C_i \}$ of curves in $M$ with 
$\lim\limits_{i \arrow \infty} \frac{\deg C_i}{2\pi(2g(C_i)-2)}\arrow \infty$. 
Suppose that $M$ is not algebraically hyperbolic. 
This gives a sequence $\phi_{C_i}:\; \Delta_1 \arrow M$ with 
$|\phi_{C_i}| > \frac{\deg C_i}{2\pi(2g(C_i)-2)}$ somewhere on $\Delta$. 
Replacing $\Delta_1$ by $\Delta_{1-\epsilon}$
as above if necessary, we may assume that $|\phi_{C_i}|$
reaches its maximum somewhere on $\Delta_1$. 
Applying an isometry of $\Delta_1$, we may assume that 
$|\phi_{C_i}|$ reaches its maximum $R_i$ in $0\in \Delta_1$. 
Rescaling $\phi_{C_i}$ by $R_i$, we obtain 
a sequence of disks 
$\tilde \phi_{C_i}(z) = \phi_{C_i}(z/R_i):\; \Delta_{R_i}\arrow M$, 
giving a Brody curve by \ref{Lemma1}. 
\endproof

\subsection{Algebraic Kobayashi metric}

\definition
The Kobayashi pseudometric $d_K$ is defined as the infimum 
of the path metric on connected chains of holomorphic 
disks equipped with the Poincar\'e metrics. 
Define the {\it algebraic Kobayashi pseudometric} $d_A$ 
using the Poincar\'e metric on connected chains 
of algebraic curves instead of the disks. 

\hfill

\remark
This metric clearly satisfies $d_A \geq d_K$. 
In \cite{Demailly_Lempert_Shiffman} Demailly, Lempert and Shiffman 
prove that $d_A = d_K$ for any quasi-projective manifold, i.e., 
one can compute the Kobayashi pseudodistance by means of chains of 
algebraic curves. 
Also, algebraic hyperbolicity follows 
if $d_A$ is non-degenerate; this is proven 
by reverse-engineering the argument we used 
to prove that algebraic hyperbolicity is implied 
by Kobayashi hyperbolicity. It is not 
clear, however, if algebraic 
hyperbolicity implies the non-degeneracy of $d_A = d_K$. 

\hfill

\question\label{_Aut_M_for_d_A>0_Question_}
Is $d_A$ always non-degenerate for algebraically
hyperbolic manifolds?

\hfill

\remark
Note that the group $\Aut(M)$ of holomorphic automorphisms
of a complex manifold $M$ with $d_A$ non-degenerate is necessarily
compact. Indeed, any holomorphic automorphism of $M$ is an isometry of
$(M, d_A)$, and the group of isometries of a compact metric space is compact.
On the other hand, a uniform limit of holomorphic maps is holomorphic,
hence compactness of the isometry group implies compactness
of $\Aut(M)$. 
Since an algebraically hyperbolic manifold has discrete 
group of holomorphic automorphisms by \ref{discrete}, the group $\Aut(M)$
is finite. Then an affirmative answer to \ref{_Aut_M_for_d_A>0_Question_}
would imply the main statement of this paper: finiteness
of $\Aut(M)$ for algebraically hyperbolic $M$.

\hfill

\remark
For any $g$ there are complex elliptic fibrations over ${\Bbb P}^1$ without 
multiple fibers where all horizontal curves (i.e., curves surjecting onto 
${\Bbb P}^1$ under the projection) have genus greater than $g$ 
(Remark 6.11 in \cite{BT-example}). The Kobayashi pseudometric on 
such surfaces is 
trivial but there are no elliptic or rational curves connecting points in 
different fibers. This example illustrates the difference
in obtaining the Kobayashi pseudometric on the surface from mapping complex
surfaces and algebraic curves. We can connect any two points in the surface
by the image of a complex line, as shown in
\cite{_Buzzard_Lu_}. However, we need to take a limit over 
all chains of algebraic curves connecting points in different fibers in 
order to obtain the Kobayashi pseudometric using only algebraic curves.

\subsection{Automorphisms of  hyperbolic manifolds}

In the rest of the section we summarize some facts about the automorphism group 
of hyperbolic manifolds. 

\hfill 

\proposition \label{discrete}
The group of automorphisms of an algebraically hyperbolic
manifold $M$ is discrete.

\hfill

{\bf Proof}. 
The group of automorphisms $G$ of a projective manifold $M$ 
is a complex Lie group. If its connected component 
$G^0$ is non-trivial, this gives a holomorphic 
map $\phi:\; G^0\arrow M$, depending on the choice of a point. 
Any connected complex algebraic Lie group is an extension 
of an affine group and an abelian variety. 
Affine algebraic groups are rational varieties, 
hence the orbit of an affine algebraic group is unirational 
and covered by rational curves. An abelian 
variety is not algebraically hyperbolic, because an 
abelian variety of dimension $n$ admits a self-isogeny of 
order $m^n$ mapping a curve $C$ of genus $g$ to a curve of genus $g$ 
and degree $m^n\deg C$. The same argument shows that any 
positive-dimensional orbit of a compact complex commutative Lie group 
is not algebraically hyperbolic. Notice that the image of an algebraically
non-hyperbolic manifold is also not algebraically hyperbolic. 
Therefore, $M$ cannot be algebraically hyperbolic. 
\endproof

\hfill

An even stronger statement is true for Kobayashi hyperbolic manifolds. 

\hfill

\proposition (Kobayashi, \cite{_Kobayashi:1976_}) 
The group of automorphisms of a compact Kobayashi hyperbolic
manifold $M$ is finite.

\hfill

{\bf Proof}. 
The group $G$ of automorphisms of $M$ 
is closed in its group of isometries under the Kobayashi pseudometric. 
The group of isometries of a compact metric space 
is compact, hence $G$ has only finitely many connected components. 
Finally, $\dim G^0=0$ as shown in \ref{discrete}, i.e., $G$ is discrete 
and compact, therefore finite. 
\endproof


\section{Main Results}


In this section $M$ is a compact complex projective manifold. 
For the classical notions of the Picard scheme and the Albanese variety 
we refer the readers to the survey manuscript \cite{Kleiman}. 

\hfill

\proposition \label{Prop1}
Let $M$ be a complex projective manifold, $\dim_\C M=n$. 
Suppose that the image of $\Aut(M)$ in $GL(H^{1,1}(M,\R))$ 
does not preserve any rational K\"ahler class. 
Then $M$ is not algebraically hyperbolic. 

\hfill 

{\bf Proof}. Assume $M$ is algebraically hyperbolic. 
Let $\omega$ be a rational K\"ahler class 
which has an infinite $\Aut(M)$-orbit. 
Replacing $\omega$ by $N\omega$, from the $(1,1)$-theorem we may 
assume that $\omega$ is a class of a hyperplane 
section. Then $\omega^{n-1} = [C]$ is the fundamental 
class of a smooth complex curve $C\subset M$ from Bertini's theorem.
Let $f_i (\omega)$ be in the 
orbit of $\omega$, which is infinite by assumption. Then 
\[ 
 \deg_\omega f_i(C)= \int_M \omega \wedge (f_i (\omega))^{n-1}=
 \int_M f_i^{-1}(\omega) \wedge \omega^{n-1}.
\]
Since the genus of $f_i(C)$ is constant, from algebraic
hyperbolicity we obtain $\int_M f_i(\omega) \wedge \omega^{n-1}<A$
for some constant $A>0$. 

Let $|\cdot |$ denote the positive-definite 
Hodge-Riemann metric on $H^{p,q}(M, \R)$, that is, the metric defined as 
$\eta, \eta' \arrow  -(-1)^l (\1)^{p-q}\int_M \eta \wedge \overline{\eta'}\wedge \omega^{n-p-q}$,
where $\eta, \eta'$ are $(p,q)$-forms which belong to a weight $l$
representation of $\goth{sl}(2)$ associated with the Lefschetz $\goth{sl}(2)$-action. 
Replace the sequence $\frac{f_i(\omega)}{|f_i(\omega)|}$ of points in the sphere 
$S\subset H^{1,1}(M, \R)$ by a converging subsequence, and let $R$ be its limit.
Since the sequence $f_i(\omega)$ is infinite, distinct and integral,
one has $\lim\limits_{i \arrow \infty} |f_i(\omega)|=\infty$. Then 
$\int_M f_i(\omega) \wedge \omega^{n-1}<A$ implies
that $\int_M R\wedge \omega^{n-1}=0$.
By the Hodge-Riemann relations, this gives 
$\int_M R\wedge R\wedge \omega^{n-2}= - |R|^2=-1$, hence
$$\lim_i \frac{\int_M f_i(\omega) \wedge f_i(\omega) \wedge 
\omega^{n-2}}{|f_i(\omega)|^2}=-1.$$ 
This is a contradiction, because $f_i(\omega)$  is K\"ahler,
hence $\int_M f_i(\omega) \wedge f_i(\omega) \wedge 
\omega^{n-2}>0$.
\endproof

\hfill

\proposition \label{Prop2}
Let $M$ be a complex projective manifold. 
Suppose that $\Aut(M)$ is infinite, but the image of $\Aut(M)$ in 
$\Aut(\Pic(M))$ is finite.  Then $M$ is not algebraically hyperbolic. 

\hfill

{\bf Proof}. From the assumption in the proposition 
we obtain that an infinite subgroup  $\Gamma\subset \Aut(M)$ 
acts trivially on $\Pic(M)$. Then it fixes a very ample 
line bundle $L \in \Pic(M)$. 
Therefore, $\Gamma$ acts on ${\Bbb P}H^0(M,L)^*$ 
preserving the image of the projective embedding 
$M \arrow {\Bbb P}H^0(M,L)^*$. 
Let $G = \bar \Gamma$ be the Zariski closure of $\Gamma$ in $PGL(H^0(M,L)^*)$. 
Since $M$ is Zariski closed and $\Gamma$-invariant, it is $G$-invariant. 
Since $\Gamma$ acts on $M$ with infinite orbits, the orbits 
of the $G$-action on $M$ are positive-dimensional. This is impossible, because 
$\Aut(M)$ is discrete as shown in \ref{discrete}. 
\endproof

\hfill

\remark \label{rem4prop3}
Notice that if the image ${\goth G}$ of $\Aut(M)$ in $GL(H^{1,1}(M,\R))$ 
does not preserve any rational K\"ahler classes, then ${\goth G}$
is infinite. Otherwise, take a finite orbit of a K\"ahler 
class and notice that its geometric center is an $\Aut(M)$-invariant 
K\"ahler class, because the convex hull of a set of K\"ahler classes 
lies inside the K\"ahler cone. 

\hfill

\remark \label{_Projective_inva_finite_action_Remark_}
When $M$ is projective and the image ${\goth G}$ of $\Aut(M)$ in $GL(H^{1,1}(M,\R))$ is finite, 
we could take a ${\goth G}$-orbit of a rational K\"ahler class,
and its geometric center is also rational. Therefore,
any projective manifold $M$ with ${\goth G}$ finite
admits an $\Aut(M)$-invariant integer K\"ahler class.

\hfill

\proposition \label{Prop3}
Let $M$ be a complex projective manifold. 
Suppose that the image of $\Aut(M)$ in $GL(H^{1,1}(M,\R))$ is finite,
but its image in $\Aut(\Pic^0(M))$ is infinite. Then $M$ is not algebraically 
hyperbolic.

\hfill

{\bf Proof}. Consider the orbit $\Aut(M) \cdot \alpha$ of a K\"ahler class 
$\alpha$. Its geometric center gives an $\Aut(M)$-invariant K\"ahler class 
$\omega$, because the convex envelope 
of a set of K\"ahler classes lies in the K\"ahler cone, as in \ref{rem4prop3}. 

The Albanese manifold $\Alb(M) = H^0(M, \Omega^1_M)^*/H_1(M,\Z)$ 
admits a natural $\Aut(M)$-invariant flat K\"ahler metric induced by the 
Hodge-Riemann form on $H^1(M, \C)$. 
Since $\Aut(M)$ acts on $\Alb(M)$ by isometries, 
the image ${\goth T}$ of $\Aut(M)$ in $\Aut(\Alb(M))$ is compact. 
The connected component of $\Aut(\Alb(M))$ 
is a subgroup $\Par(\Alb(M))=\Alb(M)\subset \Aut(\Alb(M))$
acting on the torus $\Alb(M)$ by parallel  transport.
Since ${\goth T}$ is compact, its image in the
discrete group $\frac{\Aut(\Alb(M))}{\Par(\Alb(M))}$
is finite. Let $\Gamma\subset {\goth T}$
be the kernel of the map ${\goth T}\arrow \frac{\Aut(\Alb(M))}{\Par(\Alb(M))}$.
This is a finite index subgroup in ${\goth T}$. Therefore, $\Gamma$ is infinite.

Since $\Gamma$ is infinite, its 
closure $T\subset \Par(\Alb M)$ is positive-dimensional. 
Take a smooth fiber $\Alb^{-1}(x)$ over $x\in \Alb(M)$. 
The general fiber of the real analytic map 
$\Alb^{-1}(T\cdot x)\stackrel \pi \arrow T\cdot x$ 
is smooth. Since all fibers of $\pi$ exist in dense families, 
all fibers of $\pi$ are smooth. Therefore, $\pi$ is a locally 
trivial fibration with isomorphic fibers. Consider 
the Zariski closure $T_1$ of $T$, which is a compact commutative
Lie group, that is, the compact complex torus. This gives an 
isotrivial fibration $\Alb^{-1}(X)\stackrel {\pi_1} \arrow X$, 
where $X$ is an orbit of $T_1$. 

Isotrivial fibrations over $T_1$ with fiber $F$ are classified by 
$H^1(T_1, \Aut(F))$. The variety $F\subset M$ may be assumed to be 
an algebraically hyperbolic manifold with $\dim F < \dim M$. 
Using induction on dimension, we may 
assume that $\Aut(F)$ is finite. 

The first cohomology 
of a torus with coefficients in a finite group $A$ is the same 
as an $A$-valued local system. Therefore, it becomes 
trivial after an appropriate finite covering. 
Then $\pi_1$ becomes a trivial fibration after passing to a 
finite covering $Y \arrow X$, giving 
a decomposition $\widetilde{\Alb^{-1}(X)}=F\times Y$. This manifold 
admits self-isogenies of arbitrary high order, giving 
curves of constant genus and arbitrary high degree in
$\widetilde{\Alb^{-1}(X)}$
and its finite quotient $\Alb^{-1}(X)$. 
Therefore, $\Alb^{-1}(X)\subset M $ is not algebraically hyperbolic, 
and $M$ is algebraically non-hyperbolic. 
\endproof

\hfill

Using the above results we can prove the main result of this paper. 

\hfill 

\theorem
The group $\Aut(M)$ of automorphisms of an algebraically hyperbolic
manifold $M$ is finite.

\hfill

{\bf Proof}. 
Assume $\Aut(M)$ is infinite. If the image ${\goth G}$ of the group
$\Aut(M)$ in $GL(H^{1,1}(M,\R))$ 
does not preserve any rational K\"ahler class, ${\goth G}$ is infinite
by \ref{_Projective_inva_finite_action_Remark_}. 
Then we would get a contradiction by \ref{Prop1}. 
If the image of $\Aut(M)$ in 
$\Aut(\Pic(M))$ is finite, we get a contradiction by \ref{Prop2}, and 
if its image in $\Aut(\Pic(M))$ is infinite, we get a contradiction by 
\ref{Prop3}. \endproof

\hfill

To sum the argument up: 
we consider the image ${\goth G}$ of the group
$\Aut(M)$ in $GL(H^{1,1}(M,\R))$. If it is finite,
we have an integer K\"ahler class which is ${\goth G}$-invariant.
If $\Aut(M)$ acts with finite image on the Albanese variety of $M$,
this class gives an $\Aut(M)$-invariant polarization,
hence $\Aut(M)$ is a linear algebraic group; its connected
part is trivial because the orbits of a connected 
linear algebraic group are rational. 

If ${\goth G}$ is finite, but
$\Aut(M)$ acts with infinite image on the Albanese variety of $M$,
the closure of an orbit of $\Aut(M)$ contains a torus, also
contradicting hyperbolicity. 

Finally, if ${\goth G}$ is infinite, be obtain curves of
bounded genus and arbitrary degree using the action of
${\goth G}$ on cohomology.

\hfill

{\bf Acknowledgments.} 
Most of the results in this paper were finalized during the Komplexe Analysis 
Oberwolfach Workshop in 2017. The second and third-named authors are 
grateful to the Oberwolfach organizers and staff for their hospitality. 
The first-named author has been funded by 
the Russian Academic Excellence Project '5-100' and acknowledges 
support by Simons travel grant and by the EPSRC program grant EP/M024830. 
We are thankful to the referee for many valuable comments and corrections
to the earlier version.

\noindent {\sc Fedor A. Bogomolov\\
Department of Mathematics\\
Courant Institute, NYU \\
251 Mercer Street \\
New York, NY 10012, USA,} \\
\tt bogomolov@cims.nyu.edu, also: \\
{\sc National Research University, Higher School of Economics, Moscow, Russia,}
\\

\noindent {\sc Ljudmila Kamenova\\
Department of Mathematics, 3-115 \\
Stony Brook University \\
Stony Brook, NY 11794-3651, USA,} \\
\tt kamenova@math.sunysb.edu
\\

\noindent {\sc Misha Verbitsky\\
            {\sc Instituto Nacional de Matem\'atica Pura e
              Aplicada (IMPA) \\ Estrada Dona Castorina, 110\\
Jardim Bot\^anico, CEP 22460-320\\
Rio de Janeiro, RJ - Brasil}\\
also:\\
{\sc Laboratory of Algebraic Geometry,\\
National Research University Higher School of Economics,\\
Department of Mathematics, 6 Usacheva street, Moscow, Russia.}\\
}

\end{document}